\documentstyle[11pt,amssymb]{article}

\newlength{\minitwocolumn}
\newcommand{\Z}{{\Bbb Z}} 
\newcommand{\R}{{\Bbb R}} 
\newcommand{\C}{{\Bbb C}} 

\newcommand{\F}{{\cal F}}

\newcommand{\cR}{{\cal R}}

\newcommand{\hL}{\widehat{L}}

\renewcommand{\H}{{\cal H}}

\newcommand{\nn}{{\nonumber}}
\newcommand{\eqref}[1]{(\ref{#1})}
\newcommand{\bea}{\begin{eqnarray}}
\newcommand{\ena}{\end{eqnarray}}
\newcommand{\beit}{\begin{itemize}}
\newcommand{\enit}{\end{itemize}}

\newcommand{\be}{\begin{eqnarray*}}
\newcommand{\en}{\end{eqnarray*}}
\newcommand{\lb}[1]{\label{#1}}

\def\infq4p#1{{(#1;q^4,p)_\infty}}

\newcommand{\al}{\alpha}

\newcommand{\bep}{\bar{\epsilon}}
\newcommand{\bbep}{\eta}
\newcommand{\ba}{\bar{\alpha}}

\newcommand{\bfv}{{\bf v}}

\font\teneufm=eufm10
\font\seveneufm=eufm7
\font\fiveeufm=eufm5
\newfam\eufmfam
\textfont\eufmfam=\teneufm
\scriptfont\eufmfam=\seveneufm
\scriptscriptfont\eufmfam=\fiveeufm

\let\goth\frak
\newcommand{\slth}{\widehat{\goth{sl}}_2}

\newcommand{\slnh}{\widehat{\goth{sl}}_N}
\newcommand{\sln}{\goth{sl}_N}
\newcommand{\g}{\goth{g}}

\newcommand{\Aqp}{{\cal A}_{q,p}}

\newcommand{\Bqla}{{{\cal B}_{q,\lambda}}}

\newcommand{\slnhbig}{\widehat{\mbox{\fourteeneufm sl}}_N}  

\font\fourteeneufm=eufm10 scaled\magstep2    
   
\makeatletter
\@addtoreset{equation}{section}
\makeatother

\newtheorem{thm}{Theorem}[section]
\newtheorem{prop}[thm]{Proposition}

\newtheorem{cor}[thm]{Corollary}
\newtheorem{conj}[thm]{Conjecture}

\newtheorem{df}{Definition}[section]
\newtheorem{dfn}[thm]{Definition}

\begin{document}

\

\vspace{0.7cm}
\begin{center}
{\Large{\bf
The Elliptic Algebra $U_{q,p}(\slnhbig)$ and }}

 
{\Large{\bf the Deformation of $W_N$ Algebra
}}

\vspace{1.2cm}
{ Takeo KOJIMA$~^{*}$ and~ Hitoshi KONNO$~^{**,\dagger}$}

\vspace{0.7cm}
{\it
~*
Department of Mathematics,
College of Science and Technology,\\
Nihon University, Chiyoda-ku, Tokyo
101-0062, Japan.\\
E-mail:kojima@math.cst.nihon-u.ac.jp\\~\\
~**
Department of Mathematics,
Faculty of Integrated Arts and Sciences,\\
Hiroshima University,
Higashi-Hiroshima 739-8521, Japan.
\\
E-mail:konno@mis.hiroshima-u.ac.jp
}

\vspace{0.8cm}
{\it
~$\dagger$
Department of Applied Mathematics and Theoretical Physics,\\  
Centre for Mathematical Sciences, 
University of Cambridge,\\
Wilberforce Road, Cambridge CB3 0WA, UK. }
 
\end{center}

\vspace{0.8cm}
\begin{abstract}

After reviewing the recent results on the Drinfeld realization of the 
face type elliptic quantum group $\Bqla(\slnh)$ by the elliptic algebra 
$U_{q,p}(\slnh)$, we investigate a fusion of the vertex operators of 
$U_{q,p}(\slnh)$.  The basic generating functions 
$\Lambda_j(z)\ (1\leq j\leq N-1)$ 
 of the deformed $W_N$ algebra are derived explicitly.
\end{abstract}

\newpage
\section{Introduction} 
In recent papers \cite{Konno,JKOS2,KK}, we showed that the elliptic algebra 
$U_{q,p}(\slnh)$ provides the Drinfeld realization of the face type elliptic 
quantum group 
 $\Bqla(\slnh)$ \cite{JKOS1} tensored by a Heisenberg algebra. 
Based on this fact, we defined the  $U_{q,p}(\slnh)$ counterparts of 
the intertwining operators of the 
 $\Bqla(\slnh)$ modules and 
obtained their free field realization in the level one representation. 
   The resultant vertex operators,  called the vertex operators of $U_{q,p}(\slnh)$, 
 are identified with the vertex operators of the $\slnh$ type RSOS model
  in the algebraic analysis formulation\cite{JM}. 
In general, we expect that the elliptic algebra $U_{q,p}(\g)$ with $\g$ being 
an affine Lie algebra provides the Drinfeld realization for the elliptic 
quantum group $\Bqla(\g)$ and enables us to perform 
an algebraic analysis of the $\g$ type RSOS model.   
    
On the other hand, the $\slnh$ RSOS model is known as an off-critical 
deformation of the $W_N$ minimal model\cite{JMO}. 
In this relation, it is remarkable  that the elliptic algebra 
$U_{q,p}(\slnh)$ 
in the $c=1$ representation coincides with the algebra of the screening 
currents 
of the deformed $W_N$ algebra \cite{FF,AKOS,FR}. 
In general, we expect 
that the elliptic algebra $U_{q,p}(\g)$ provides an algebra of screening 
currents of the deformation of the coset CFT  
associated with  $(\g)_c\oplus (\g)_{r-c-2}/(\g)_{r-2}$\cite{Konno,JKOS2}, 
which corresponds to the $c\times c$ fusion RSOS model of type $\g$.

The purpose of this paper is to continue to 
discuss an explicit relation among the elliptic algebra $U_{q,p}(\g)$,
 the $\g$ type RSOS model 
and the deformation of $W(\bar{\g})$ algebra in the case $\g=\slnh$.
We here investigate a fusion of the type II vertex operator of $U_{q,p}(\slnh)$  and its  dual, and   
show that the generating functions of the deformed $W_N$ 
algebra can be extracted from it. 
The idea of fusion of the vertex operators was used in \cite{JKM,JS}  
to  derive the generating function of the deformed Virasoro algebra ( corresponding to the $\phi_{1,3}$ perturbation) from the ABF model in regime III, in 
\cite{JKOPS}for the deformed $W_N$ algebra with the central charge $c_N=(N-1)\left(1-\frac{N(N+1)}{r(r-1)}\right)$ at special point 
$r=N+2$ (the $\Z_N$ parafermion point) from the ABF model in regime II, and in \cite{HJKOS} for the deformed Virasoro algebra ( corresponding to the $\phi_{1,2}$ perturbation) from the dilute $A_L$ model. 

This article is organised as follows. In the next section, we briefly review
 the elliptic algebra $U_{q,p}(\slnh)$ as the Drinfeld realization of the 
 elliptic quantum group $\Bqla(\slnh)$ according to \cite{KK}. In section 3, we give a summary of the results on   the free field realization of 
 the vertex operators of $U_{q,p}(\slnh)$. 
In section 4, we discuss a fusion of the vertex operators 
of $U_{q,p}(\slnh)$ and derive the basic generators of the deformed $W_N$ 
algebra.

Through this paper, we use  the  following symbols.
$p=q^{2r},~~p^*=pq^{-2c}=q^{2r^*}~~(r^*=r-c;~ r,r^* \in 
{\mathbb{R}}_{>0})$, 
\begin{eqnarray*}
&&[n]_q=\frac{q^n-q^{-n}}{q-q^{-1}},\\
&&\Theta_p(z)=(z,p)_\infty (pz^{-1};p)_\infty
(p;p)_\infty,\\
&&\{z\}=(z;p,q^{2N})_{\infty},\qquad \{z\}^*=\{z\}|_{p\to p^*},\\
&&(z;t_1,\cdots,t_k)_\infty=
\prod_{n_1,\cdots,n_k \geqq 0}(1-zt_1^{n_1}\cdots t_k^{n_k}).
\end{eqnarray*}
We also use the Jacobi theta functions
\begin{eqnarray*}
[v]=q^{\frac{v^2}{r}-v}
\frac{\Theta_p(q^{2v})}{(p;p)_\infty^3},~~
[v]^*=q^{\frac{v^2}{r^*}-v}
\frac{\Theta_{p^*}(q^{2v})}{(p^*;p^*)_\infty^3},
\end{eqnarray*}
which satisfy $[-v]=-[v]$ and the quasi-periodicity property
\begin{eqnarray*}
&&~[v+r]=-[v],~~[v+r\tau]=-e^{-\pi i \tau-\frac{2\pi i v}{r}}[v].
\end{eqnarray*}
We take the normalization of the theta function to be 
\begin{eqnarray*}
\oint_{C_0}\frac{dz}{2\pi i z}\frac{1}{[-v]}=1,\qquad \oint_{C_0}\frac{dz}{2\pi i z}\frac{1}{[-v]^*}=\frac{[v]}{[v]^*}
\Bigl\vert_{v\to 0}.
\end{eqnarray*}
where $C_0$ is a simple closed curve in the $v$-plane encircling
$v=0$ anticlockwise.

\section{The Elliptic Algebra $U_{q,p}(\slnhbig)$}\lb{secuqpslnh}

\subsection{Definition}\lb{drinfelduqp}

\begin{dfn}{\bf ( Elliptic algebra $U_{q,p}(\slnh)$ )}\lb{uqpslnh}
We define the elliptic algebra $U_{q,p}(\slnh)$ to be the 
associative algebra of the currents $E_j(v),\ F_j(v)\ (1\leq j\leq N-1)$ and 
$ K_j(v)\ (1\leq j\leq N)$ satisfying the following relations.
\end{dfn}
\begin{eqnarray}
&&E_i(v_1)E_j(v_2)=
\frac{[v_1-v_2+\frac{A_{ij}}{2}]^*}{[v_1-v_2-\frac{A_{ij}}{2}]^*}
E_j(v_2)E_i(v_1),
\label{EiEj}\\
&&\nn\\
&&F_i(v_1)F_j(v_2)=\frac{[v_1-v_2-\frac{A_{ij}}{2}]}
{[v_1-v_2+\frac{A_{ij}}{2}]}
F_j(v_2)F_i(v_1),
\label{FiFj}\\
&&\nn\\
&&[E_i(v_1),F_j(v_2)]=\frac{\delta_{i,j}}{q-q^{-1}}
\left(\delta(q^{-c}z_1/z_2)H_j^+\left(v_2+\frac{c}{4}\right)
-\delta(q^{c}z_1/z_2)H_j^-\left(v_2-\frac{c}{4}\right)
\right),
\label{EiFj}\\
&&\nn\\
&&H_j^\pm\left(v\mp \frac{1}{2}\left(r-\frac{c}{2}\right)\right)
=\kappa K_j\left(v+\frac{N-j}{2}\right)
K_{j+1}\left(v+\frac{N-j}{2}\right)^{-1},\label{HKK}\\
&&\nn\\
&&{K_j}(v_1){K_j}(v_2)={\rho^{}(v_1-v_2)}{K_j}(v_2){K_j}(v_1)
,\label{k2}
\\
&&{K_{j_1}}(v_1){K_{j_2}}(v_2)=
{\rho(v_1-v_2)}
\frac{[v_1-v_2-1]^*[v_1-v_2]}{[v_1-v_2]^*[v_1-v_2-1]}
{K_{j_2}}(v_2){K_{j_1}}(v_1)
\nonumber\\
&&~~~~~~~~~~~~~~~~\qquad\qquad (1\leqq j_1<j_2\leqq N),\label{k3}\\
&&{K_j}(v_1){E_j}(v_2)=
\frac{[v_1-v_2+\frac{j+r^*-N}{2}]^*}
{[v_1-v_2+\frac{j+r^*-N}{2}-1]^*}
{E_j}(v_2){K_j}(v_1),\label{k4}\\
&&{K_{j+1}}(v_1){E_j}(v_2)=
\frac{[v_1-v_2+\frac{j+r^*-N}{2}]^*}
{[v_1-v_2+\frac{j+r^*-N}{2}+1]^*}
{E_j}(v_2){K_{j+1}}(v_1),\label{k5}\\
&&{K_{j_1}}(v_1){E_{j_2}}(v_2)=
{E_{j_2}}(v_2){K_{j_1}}(v_1)~~\qquad (j_1\neq j_2,j_2+1),
\label{k6}\\
&&\nn\\
&&{K_j}(v_1){F_j}(v_2)=
\frac{[v_1-v_2+\frac{j+r-N}{2}-1]}
{[v_1-v_2+\frac{j+r-N}{2}]}
{F_j}(v_2){K_j}(v_1),
\label{k7}
\\
&&K_{j+1}(v_1)F_j(v_2)=
\frac{[v_1-v_2+\frac{j+r-N}{2}+1]}
{[v_1-v_2+\frac{j+r-N}{2}]}F_j(v_2)K_{j+1}(v_1)
,\label{k8}
\\
&&K_{j_1}(v_1)F_{j_2}(v_2)=
F_{j_2}(v_2)K_{j_1}(v_1)~~\qquad (j_1\neq j_2,j_2+1),\label{k9}\\
&&\nn\\
%
&&z_1^{-\frac{1}{r^*}}
\frac{(p^*q^2z_2/z_1;p^*)_\infty}
{(p^*q^{-2}z_2/z_1;p^*)_\infty
}\left\{
({z_2}/{z})^{\frac{1}{r^*}}
\frac{(p^*q^{-1}z/z_1;p^*)_\infty 
(p^*q^{-1}z/z_2;p^*)_\infty}
{(p^*qz/z_1;p^*)_\infty 
(p^*qz/z_2;p^*)_\infty}E_i(v_1)E_i(v_2)E_j(v)\right.\nonumber\\
&&
-\left.[2]_q\frac{(p^*q^{-1}z/z_1;p^*)_\infty 
(p^*q^{-1}z_2/z;p^*)_\infty}
{(p^*qz/z_1;p^*)_\infty 
(p^*qz_2/z;p^*)_\infty}
E_i(v_{1})E_j(v)E_i(v_{2})
\right.\nn
\\
&&+\left.
(z/z_1)^{\frac{1}{r^*}}\frac{(p^*q^{-1}z_1/z;p^*)_\infty 
(p^*q^{-1}z_2/z;p^*)_\infty}
{(p^*qz_1/z;p^*)_\infty 
(p^*qz_2/z;p^*)_\infty}
E_j(v)E_i(v_{1})E_i(v_{2})
\right\}+(z_1 \leftrightarrow z_2)=0,\nn\\
&&\label{e20}\\
&&z_1^{\frac{1}{r}}
\frac{(pq^{-2}z_2/z_1;p)_\infty}
{(pq^{2}z_2/z_1;p)_\infty
}\left\{
(z/z_2)^{\frac{1}{r}}
\frac{(pq z/z_1;p)_\infty 
(pq z/z_2;p)_\infty}
{(pq^{-1} z/z_1;p)_\infty 
(pq^{-1} z/z_2;p)_\infty}
F_i(v_1)F_i(v_2)F_j(v)\right.\nonumber\\
&&
-\left.[2]_q\frac{(pq z/z_1;p)_\infty 
(pq z_2/z;p)_\infty}
{(pq^{-1} z/z_1;p)_\infty 
(pq^{-1} z_2/z;p)_\infty}
F_i(v_{1})F_j(v)F_i(v_{2})
\right.\nn
\\
&&+\left.
(z_1/z)^{\frac{1}{r}}
\frac{(pq z_1/z;p)_\infty 
(pq z_2/z;p)_\infty}
{(pq^{-1} z_1/z;p)_\infty 
(pq^{-1} z_2/z;p)_\infty}
F_j(v)F_i(v_{1})
F_i(v_{2})
\right\} + (z_1 \leftrightarrow z_2)=0
\quad (|i-j|=1).\nonumber\\
&&~~~~~~~~~~~~~~~~~\label{e21}
\end{eqnarray}
Here $A=(A_{j k})$ is the Cartan matrix of $\sln$. The constant $\kappa$ and  
the functions $\rho(v)$ are given by 
\begin{eqnarray}
&&{\kappa}=\frac{(p;p)_\infty (p^*q^2;p^*)_\infty}
{(p^*;p^*)_\infty (pq^2;p)_\infty},\\
&&\rho(v)=\frac{\rho^{+*}(v)}{\rho^+(v)},\\
&&\rho^+(v)=q^{\frac{N-1}{N}}
z^{\frac{N-1}{rN}}
\frac{\{pq^2z\}
\{pq^{2N-2}z\}
\{1/z\} \{q^{2N}/z\}
}{
\{pz\}
\{pq^{2N}z\}
\{q^2/z\} \{q^{2N-2}/z\}
},\quad 
\rho^{+*}(v)=\rho^+(v)|_{r\to r^*}.
\label{def:rho*}
\end{eqnarray}

\subsection{Realization of $U_{q,p}(\slnh)$}
The elliptic algebra $U_{q,p}(\slnh)$ can be realized by using the Drinfeld 
generators of $U_q(\slnh)$ and a Heisenberg algebra. 
Let $h_{i},\  a^i_{m},\ x_{i,n}^\pm$ $(1\leq i \leq N-1: m 
\in {{\Z}_{\neq 0}},\ n \in {\mathbb{Z}}),\ c,\ d$ be the standard Drinfeld generators of $U_q(\slnh)$\cite{Drinfeld,KK}. Their generating functions $x_i^\pm(z)$, $\psi_i(z)$, $\varphi_i(z)$ are called the Drinfeld currents. 
\begin{eqnarray}
&&x_i^\pm(z)=\sum_{n \in {\mathbb{Z}}}x_{i,n}^\pm z^{-n},\\
&&\psi_i(q^{\frac{c}{2}}z)=q^{h_i}
\exp\left(
(q-q^{-1})\sum_{m>0}a_{i,m}z^{-m}\right),\\
&&\varphi_i(q^{-\frac{c}{2}}z)=q^{-h_i}
\exp\left(-
(q-q^{-1})\sum_{m>0}a_{i,-m}z^{m}\right)\qquad (1\leq i\leq N-1). 
\end{eqnarray}

\begin{df}
We define ``dressed'' currents
${e}_i(z,p), {f}_i(z,p)$, 
${\psi}_i^\pm(z,p),$
$(1\leq i\leq N-1)$ by
\begin{eqnarray}
&&e_i(z,p)=u_i^+(z,p)x_i^+(z),\label{def:e}\\
&&f_i(z,p)=x_i^-(z)u_i^-(z,p),\label{def:f}\\
&&\psi_i^+(z,p)=u_i^+(q^{\frac{c}{2}}z,p)\psi_i(z)
u_i^-(q^{-\frac{c}{2}}z,p)
,\label{def:psi+}\\
&&\psi_i^-(z,p)=
u_i^+(q^{-\frac{c}{2}}z,p)\varphi_i(z)
u_i^-(q^{\frac{c}{2}}z,p),\label{def:psi-}
\end{eqnarray}
where 
\begin{eqnarray}
&&u_i^+(z,p)=\exp\left(
\sum_{m>0}\frac{1}{[r^* m]_q}a_{i,-m}(q^r z)^m\right),
\label{def:u1}\\
&&u_i^-(z,p)=
\exp\left(
-\sum_{m>0}\frac{1}{[r m]_q}a_{i,m}(q^{-r}z)^{-m}\right).
\label{def:u2}
\end{eqnarray}
\end{df}

Setting $\displaystyle{
b_{j,m}=
\frac{[r^*m]_q}{[rm]_q}a_{j,m}}$ (for $ m>0$), 
$\displaystyle{q^{c|m|}a_{j,m}}$ (for $m<0$), 
we introduce new generators,
$B_m^j\ (1\leq j\leq N; m \in {\mathbb{Z}})$, by 
\begin{eqnarray}
-B_m^j+B_m^{j+1}=\frac{m}{[m]_q}b_{j,m} q^{(N-j)m},~~\qquad
\sum_{j=1}^N q^{2jm} B_m^j=0.\label{def:Boson1}
\end{eqnarray}
From this and the commutation relation of the Drinfeld generators $a_{j,m}$, 
we derive the following commutation relations.
\begin{eqnarray}
&&~[B_m^j,B_{m'}^k]=m \delta_{m+m',0}
\frac{[r^*m]_q[cm]_q}{[rm]_q[m]_q[Nm]_q}\times
\left\{\begin{array}{cc}
[(N-1)m]_q &\ (j=k)\\
-q^{-mN{\rm sgn(j-k)}}[m]_q&\ (j\neq k)
\end{array}\right.,
\end{eqnarray}
for $m,m' \in {{\Z}_{\neq 0}},\ 1\leq j,k \leq N$.
Then defining new currents $k_j(z,p)\ (1\leq j\leq N)$  by
\begin{eqnarray}
k_j(z,p)= 
:\exp\left(\sum_{m\neq0}\frac{[m]_q}
{m [r^*m]_q}B_m^j z^{-m}\right):, 
\label{def:k}
\end{eqnarray}
we obtain the following decomposition.
\begin{eqnarray}
&&\psi^\pm_j(q^{\pm (r - \frac{c}{2})}z,p)=
{\kappa}\ q^{\pm h_j} k_j(q^{N-j}z,p) 
k_{j+1}(q^{N-j}z,p)^{-1}.\lb{psikk}
\end{eqnarray} 

On the other hand, let $ \epsilon_j\ (1\leq j\leq N)$
be the orthonormal basis in $\R^N$ with the inner product
$\langle \epsilon_j, \epsilon_k \rangle=\delta_{j,k}$.
Setting $\displaystyle{\bar{\epsilon}_j=\epsilon_j-\epsilon,~\ 
\epsilon=\frac{1}{N}\sum_{j=1}^N \epsilon_j,}$ 
we have the weight lattice ${P}$ of type $A_{N-1}^{(1)}$;\    
$\displaystyle{P=\oplus_{j=1}^N {\mathbb{Z}}\ \bar{\epsilon}_j.}$
Then, for example, the simple roots  $\alpha_j\ (1\leq j\leq N-1)$ of $\sln$ are given by $\displaystyle{ 
\alpha_j=-\bar{\epsilon}_j+
\bar{\epsilon}_{j+1}.}$ 
Let us introduce operators $h_\alpha, \beta\ (\alpha,\beta \in {P})$
by 
\bea
&&~[h_{\bar{\epsilon}_j},\bar{\epsilon}_k]
=\langle \bar{\epsilon}_j,
\bar{\epsilon}_k \rangle,\qquad [h_{\bar{\epsilon}_j},h_{\bar{\epsilon}_k}]=0=
[\bar{\epsilon}_j,\bar{\epsilon}_k],\lb{HA1}
\ena
$h_\alpha=\sum_j n_j h_{\bep_j}$ for $\alpha=\sum_j n_j \bep_j$ and $h_0=0$.
Note that $\langle \bar{\epsilon}_j, \bar{\epsilon}_k \rangle=
\delta_{j,k}-\frac{1}{N}$ and $[h_{\alpha_j}, \alpha_k]=2\delta_{j,k}-\delta_{j,k+1}-\delta_{j,k-1}=A_{j k}$. 
We hence identify $h_{\alpha_j}=-h_{\bar{\epsilon}_j}+
h_{\bar{\epsilon}_{j+1}}$ with $h_j$ 
in the Drinfeld generators of $U_{q}(\slnh)$ 
\begin{dfn}
We define the ( centrally extended ) Heisenberg algebra
$\C\{\hat{\H}\}$ as an associative algebra generated by 
$P_{\bep_j},\ Q_{\bep_j}\ (1\leq j\leq N)$ and  $\bbep_j\ (1\leq j\leq N-1)$
with the relations 
\bea
&&[P_{\bar{\epsilon}_j}, Q_{\bar{\epsilon}_k}]=
\langle \bar{\epsilon}_j, \bar{\epsilon}_k \rangle, \qquad 
[P_{\bar{\epsilon}_j}, P_{\bar{\epsilon}_k}]=0,
\lb{central1}
\\
&&
[Q_{\bar{\epsilon}_{j}}, Q_{\bar{\epsilon}_{k}}]=
\left(\frac{1}{r}-\frac{1}{r^*}\right){\rm sgn}(j-k)
{\rm log}~q,\lb{central2}\\
&&[Q_{\bar{\epsilon}_{j}},\bbep_k]=
\frac{1}{r}{\rm sgn}(j-k)\log~q,\lb{central3}\\
&&[\bbep_j,\bbep_k]=\frac{1}{r}{\rm sgn}(j-k)\log~q,\lb{central4}\\
&&[P_{\bar{\epsilon}_j}, \bbep_k]=0,\qquad \sum_{j=1}^N \bbep_j=0,\lb{central5}
\\
&&
[\bbep_j, \al]=
[P_{\bar{\epsilon}_j}, U_q(\widehat{\goth{sl}}_N)]=
[Q_{\bar{\epsilon}_j}, U_q(\widehat{\goth{sl}}_N)]=
[\bbep_j, U_{q}(\slnh)]=0.
\lb{central6}
\ena
\end{dfn}

\begin{dfn}
We  define the currents 
$E_j(v), F_j(v), H_j^\pm(v)\ (1\leq j\leq N-1)$ and 
$K_j(v)\ (1\leq j\leq N)$ by 
\begin{eqnarray}
&&{E}_j(v)=e_j(z,p)e^{\bar{\alpha}_j}e^{-Q_{\alpha_j}}(q^{-j+N}z)
^{-\frac{P_{\alpha_j}-1}{r^*}},\label{def:total1c}\\
&&{F}_j(v)=f_j(z,p)e^{-\bar{\alpha}_j}
(q^{-j+N}z)^{\frac{P_{\alpha_j}-1}{r}}
(q^{-j+N}z)^{\frac{h_j}{r}},\label{def:total2c}\\
&&{H}_j^\pm
\left(v\right)=
\psi_j^\pm(z,p)
q^{\mp h_j}e^{-Q_{\alpha_j}}
(q^{-j+N\pm(r-\frac{c}{2})}z)^{(-\frac{1}{r^*}+\frac{1}{r})
(P_{\alpha_j}-1)+\frac{1}{r}h_j},\label{def:total3c}
\\
&&{K}_j(v)=k_j(z,p)e^{Q_{\bar{\epsilon}_j}}
z^{(\frac{1}{r^*}-\frac{1}{r})P_{\bar{\epsilon}_j}}
z^{-\frac{1}{r}h_{\bar{\epsilon}_j}+(\frac{1}{r^*}-\frac{1}{r})\frac{N-1}{2N}},\label{def:total4c}
\end{eqnarray}
where $z=q^{2v}$ and $\ba_j=-\bbep_j+\bbep_{j+1}$. 
\end{dfn}

Then it is easy to show that 
$E_j(v), F_j(v), H_j^\pm(v)$ and $K_j(v)$ satisfy the defining relations of 
the elliptic algebra $U_{q,p}(\slnh)$.

\subsection{$RLL$ relation}\lb{halfcurrents}
We next discuss a relation between two elliptic algebras   
$U_{q,p}(\slnh)$ and ${\cal B}_{q,\lambda}(\slnh)$.  
We construct a $L$-operator 
 by using the half currents
and show that it satisfies the dynamical $RLL$-relation 
which characterizes ${\cal B}_{q,\lambda}(\slnh)$.
We use the following abbreviations
\bea
P_{j,l}&=&
-P_{\bar{\epsilon}_{j}}+P_{\bar{\epsilon}_{l}}
=P_{\alpha_{j}}+P_{\alpha_{j+1}}+\cdots
+P_{
\alpha_{l-1}}~~\\
h_{j,l}&=&-h_{\bar{\epsilon}_j}+h_{\bar{\epsilon}_l}=h_{j}+h_{j+1}+\cdots+
h_{l-1}~~
\end{eqnarray}
for $ j<l $. From the definition of $\C\{\hat{\H}\}$ and \eqref{def:total1c}-\eqref{def:total4c}, we have
\begin{eqnarray}
&&~[K_j(v),P_{k,l}]=(\delta_{j,k}-\delta_{j,l})
K_j(v)=[K_j(v),P_{k,l}+h_{k,l}],\lb{kjpjl}\\
&&~[E_j(v),P_{k,l}]=(\delta_{j,k}+
\delta_{j+1,l}-\delta_{j,l}-\delta_{j+1,k})E_j(v),\lb{ejpjl}\\
&&~[F_j(v),P_{j,l}+h_{j,l}]=(
\delta_{j,k}+
\delta_{j+1,l}-\delta_{j,l}-\delta_{j+1,k}
)F_j(v),\lb{fjpjlphjl}\\
&&
~[F_j(v),P_{k,l}]=0=
[E_j(v),P_{k,l}+h_{k,l}].\lb{kakan}
\end{eqnarray}

\begin{df}
We define the half currents
$F_{j,l}^+(v), E_{l,j}^+(v), (1\leq j<l \leq N)$
and $K_j^+(v)\ (1\leq j \leq N)$
by
\begin{eqnarray}
K_j^+(v)&=&K_j\left(v+\frac{r+1}{2}\right)~~~
(1\leq j \leq N),\lb{kjp}\\
F_{j,l}^+(v)&=&a_{j,l}\oint_{C(j,l)}
\prod_{m=j}^{l-1}\frac{dw_m}{2\pi i w_m}
F_{l-1}(v_{l-1})F_{l-2}(v_{l-2})
\cdots F_{j}(v_{j})\nonumber\\
&&\qquad\times
\frac{[v-v_{l-1}+P_{j,l}+h_{j,l}+\frac{l-N}{2}-1]
[1]}{[v-v_{l-1}+\frac{l-N}{2}]
[P_{j,l}+h_{j,l}-1]}\nonumber\\
&&\qquad\times
\prod_{m=j}^{l-2}
\frac{[v_{m+1}-v_{m}+P_{j,m+1}+h_{j,m+1}-\frac{1}{2}][1]}
{[v_{m+1}-v_{m}+\frac{1}{2}][P_{j,m+1}+h_{j,m+1}
]},\lb{fjl}\\
E_{l,j}^+(v)&=&a_{j,l}^*\oint_{C^*(j,l)}
\prod_{m=j}^{l-1}\frac{dw_m}{2\pi i w_m}
E_{j}(v_{j})E_{j+1}(v_{j+1})\cdots E_{l-1}(v_{l-1})
\nonumber\\
&&\qquad\times
\frac{[v-v_{l-1}-P_{j,l}+\frac{l-N}{2}+
\frac{c}{2}+1]^*[1]^*}
{[v-v_{l-1}+\frac{l-N}{2}+\frac{c}{2}]^*[P_{j,l}-1]^*
}
\nonumber\\
&&\qquad\times
\prod_{m=j}^{l-2}
\frac{[v_{m+1}-v_{m}-P_{j,m+1}+\frac{1}{2}]^*[1]^*}
{[v_{m+1}-v_{m}+\frac{1}{2}]^*[P_{j,m+1}-1]^*}.\lb{elj}
\end{eqnarray}  
Here $w_m=q^{2v_m}$ and the integration contour $C(j,l)$ and $C^*(j,l)$ are given by
\bea
C(j,l)&:& |pq^{l-N}z|<|w_{l-1}|<|q^{l-N}z|,\nn\\
&& |pqw_{m+1}|<|w_{m}|<|qw_{m+1}|,\\
C^*(j,l)&:& |p^*q^{l-N+c}z|<|w_{l-1}|<|q^{l-N+c}z|,\nn\\
&& |p^*qw_{m+1}|<|w_{m}|<|qw_{m+1}|,
\ena
where $m=j,j+1,..,l-2$.
The constants $a_{j,l}$ and $a_{j,l}^*$ are chosen to 
satisfy
\begin{eqnarray}
&&\frac{\kappa\ a_{j,l}a_{j,l}^* [1]}{q-q^{-1}}=1.
\end{eqnarray}
\end{df}


\subsection{$L$-operator}
\begin{df}
By using the half currents, 
we define the $L$-operator
$\widehat{L}^+(v) \in {\rm End}({\mathbb{C}}^N)
\otimes U_{q,p}(\slnh)$
as follows.
\begin{eqnarray}
&&\widehat{L}^+(z)=
\left(\begin{array}{ccccc}
1&F_{1,2}^+(v)&F_{1,3}^+(v)&\cdots&F_{1,N}^+(v)\\
0&1&F_{2,3}^+(v)&\cdots&F_{2,N}^+(v)\\
\vdots&\ddots&\ddots&\ddots&\vdots\\
\vdots&&\ddots&1&F_{N-1,N}^+(v)\\
0&\cdots&\cdots&0&1
\end{array}\right)\left(
\begin{array}{cccc}
K^+_1(v)&0&\cdots&0\\
0&K^+_2(v)&&\vdots\\
\vdots&&\ddots&0\\
0&\cdots&0&K^+_{N}(v)
\end{array}
\right)\nn\\
&&\qquad\qquad\qquad\qquad\qquad\qquad\qquad\times
\left(
\begin{array}{ccccc}
1&0&\cdots&\cdots&0\\
E^+_{2,1}(v)&1&\ddots&&\vdots\\
E^+_{3,1}(v)&
E^+_{3,2}(v)&\ddots&\ddots&\vdots\\
\vdots&\vdots&\ddots&1&0\\
E^+_{N,1}(v)&E^+_{N,2}(v)
&\cdots&E^+_{N,N-1}(v)&1
\end{array}
\right).\lb{def:lhat}
\end{eqnarray}
\end{df}
Then a direct comparison with the relations 
of the half currents leads us to  the following conjecture.

\noindent
\begin{conj}\lb{main}~~
The $L$-operator $\widehat{L}^+(v)$ satisfies
the following $RLL=LLR^*$ relation.
\begin{eqnarray}
R^{+(12)}(u_1-u_2,P+h)\widehat{L}^{+(1)}(z_1)
\widehat{L}^{+(2)}(z_2)=
\widehat{L}^{+(2)}(z_2)
\widehat{L}^{+(1)}(z_1)
R^{+*(12)}(u_1-u_2,P).\label{thm:RLL}
\end{eqnarray}
Here $z_i=q^{2u_i}\ (i=1,2)$. The $R$-matrix $R^+(v,P)$ is the 
image of the universal $R$-matrix $\cR(r,\{s_j\})$ of $\Bqla(\slnh)$ 
in the evaluation representation $(\pi_{V_z}\otimes \pi_{V_1})$, 
$V\cong \C^N$, given by  
\begin{eqnarray}
R^{+}(v,s)&=&\rho^{+}(v)\bar{R}(v,s),\lb{rmatfull}\\
\bar{R}(v,s)&=&\sum_{j=1}^{N}
E_{jj}\otimes E_{jj}+\sum_{1 \leq j<l \leq N}
\left(b_{}(v,s_{j,l})E_{jj}
\otimes E_{ll}+\bar{b}_{}(v)E_{ll}\otimes E_{jj}\right)\nonumber\\
&&+\sum_{1 \leq j<l \leq N}\left(c_{}(v,s_{j,l})E_{jl}\otimes E_{lj}+
\bar{c}_{}(v,s_{j,l})E_{lj}\otimes E_{jl}\right),\lb{rmat}
\end{eqnarray}
where $s_{j,l}=\sum_{m=j}^{l-1}s_j \ (1\leq j<l \leq N)$ and 
\begin{eqnarray}
b_{}(u,s)&=&\frac{[s+1][s-1][u]}{[s]^{2} [u+1]},
~~~\bar{b}_{}(u)=\frac{[u]}{[u+1]},\\
c_{}(u,s)&=&\frac{[1][s+u]}{[s][u+1]},
~~~\bar{c}_{}(u,s)=\frac{[1][s-u]}{
[s][u+1]}.\lb{rmatcomp}
\end{eqnarray}
And $R^{+*}(v,s)=R^{+}(v,s)|_{r\to r^*}$. 
Up to a gauge transformation, $R^+(v,P)$ coincides with the 
Boltzmann weight of the $\slnh$ RSOS model\cite{JMO}.
\end{conj}
The  $c=1$ case, the statement was proved by using the free field realization
\cite{KK}.


Now let us define the modified $L$-operator $L^+(v,P)$  by
\begin{eqnarray}
L+(z,P)=
\widehat{L}+(z)\left(\begin{array}{cccc}
e^{-Q_{\bar{\epsilon}_1}}&0&\cdots&0\\
0&e^{-Q_{\bar{\epsilon}_2}}&&\vdots\\
\vdots&&\ddots&0\\
0&\cdots&0&e^{-Q_{\bar{\epsilon}_N}}
\end{array}\right)=\widehat{L}^+(z)\ 
\exp\left\{\sum_{m=1}^N h^{(1)}_{{\epsilon}_m}Q_{\bar{\epsilon}_m}\right\}.
\lb{lhat}
\end{eqnarray}
Here $h^{(1)}_{{\epsilon}_j}=h_{{\epsilon}_j}\otimes 1,\ h_{\epsilon_m}\equiv
-E_{mm}$ (a $N \times N$ { matrix unit}). 
We then show that the modified $L$-operator depends on neither 
$Q_{\bar{\epsilon}_j}$ nor $ {\bbep}_j$ and satisfies
the dynamical $RLL$ relation of $\Bqla(\slnh)$\cite{JKOS1}.
\begin{cor} 
\begin{eqnarray}
&&R^{+(12)}(v,P+h)L^{+(1)}(z_1,P)L^{+(2)}(z_2,P+h^{(1)})
=
L^{+(2)}(z_2,P)L^{+(1)}(z_1,P+h^{(2)})
R^{+*(12)}(v,P),\nn\\
&&\label{thm:dRLL}
\end{eqnarray}
where $u=u_1-u_2$.
\end{cor}
Hence, we regard the elliptic currents 
$E_j(v),\ F_j(v)\ (1\leq j\leq N-1)$ and 
$K_j(v)\ (1\leq j\leq N)$ in $U_{q,p}(\slnh)$ 
as the Drinfeld realization of the elliptic 
algebra $\Bqla(\slnh)$ tensored by the Heisenberg algebra. 
\begin{eqnarray}
U_{q,p}(\slnh)=\Bqla(\widehat{\goth{sl}}_N)
\otimes_{{\mathbb{C}}\{P_{\bep_1},P_{\bep_2},..,P_{\bep_{N-1}}\}} 
{\mathbb{C}}\{\hat{\H}\}.
\end{eqnarray}


\section{Vertex Operators of $U_{q,p}(\slnhbig)$}
\lb{vertexoperators}
We here summarize a construction of the type II vertex operator of $U_{q,p}(\slnh)$ 
and its dual vertex operator.

\subsection{Definition}
Let us first define an extension of the $U_q(\slnh)(\cong \Bqla(\slnh))$ modules by
\be
&&\widehat{\cal F}=\bigoplus_{\mu_1,\cdots,\mu_{N-1} \in {\mathbb{Z}}}
{\cal F}\otimes e^{\mu_1 Q_{\bar{\epsilon}_1}+
\cdots +\mu_{N-1} Q_{\bar{\epsilon}_{N-1}}}.
\en 
Let 
${\Psi}^*_W(z,P)$ be the 
type II intertwining operator of $\Bqla(\slnh)$ \cite{JKOS1}.
We define the 
type II vertex operator
 $\widehat{\Psi}^*_W(z)$ of $U_{q,p}(\slnh)$
as the following extension.
\begin{eqnarray}
&&\widehat{\Psi}^*_W(z)={\Psi}^*_W(z,P)
\exp\left\{ \sum_{j=1}^N h^{(1)}_{{\epsilon}_j} Q_{\bar{\epsilon}_j}\right\}\qquad :
W_z \otimes \widehat{\cal F} \longrightarrow \widehat{\cal F}'.
\end{eqnarray}
From the intertwining relation of the $\Bqla(\slnh)$ intertwining 
operators, we derive the following relation for 
the new operator 
$\widehat{\Psi}_W^*(z)$. 
\begin{eqnarray}
\widehat{L}_V^{+(1)}(z_1)\widehat{\Psi}_W^{*(2)}(z_2)&=&
\widehat{\Psi}_W^{*(2)}(z_2)\widehat{L}_V^{+(1)}(z_1)
R_{VW}^{+*(12)}(u_1-u_2,P-h^{(1)}-h^{(2)}).\label{Inter:Uqp2}
\end{eqnarray}

Let us consider the vector representation $V=W\cong \C^N$ of $\Bqla(\slnh)$.
We denote  a basis of $V$ by $\{\bfv_m \}_{m=1}^N$.
In this representation, the $R$-matrix $R^+_{VV}(v,P)$ is given by 
$R^+(v,P)$ in (\ref{rmatfull}) and the $L$-operator 
$\widehat{L}^+_V(z)$ by $\hL^+(z)$ in (\ref{def:lhat}).
We define the components of the vertex operators by
\begin{eqnarray}
\widehat{\Psi}_V^*\left(q^{-c-1}z\right)(\bfv_m \otimes \cdot)=
\Psi^*_m(z),
\end{eqnarray}
and the matrix elements of
the $L$-operator $\widehat{L}^+(z)$ by
\begin{eqnarray}
\widehat{L}^+(z)\bfv_{j}=
\sum_{1\leq m \leq N}
\bfv_m L^+(z)_{mj}.
\end{eqnarray}

\subsection{Free field realizations}
We here construct a free
field realization of the vertex operators 
fixing $c=1$. 
Let $\al_j$ be the simple root operator. We make the
 standard central extension $
[\al_j, \al_k]=\pi i A_{jk}$ and set $\hat{\al}_j=\al_j+\ba_j$,
 where $\ba_j$ is an element of 
the Heisenberg algebra $\C\{\hat{H}\}$. 
Then we have 
\begin{prop}~~
The currents $E_j(v)$ and $F_j(v)$ given by
\begin{eqnarray}
E_j(v)&=&:\exp\left(
-\sum_{m \neq 0}\frac{[r m]_q}{m[r^* m]_q}(-B_m^j+B_m^{j+1})
(q^{N-j}z)^{-m}\right):
e^{\hat{\alpha}_j}z^{h_j}e^{-Q_{\alpha_j}}(q^{-j+N}z)
^{-\frac{P_{\alpha_j}-1}{r^*}},\nonumber\\
\label{free:E}\\
F_j(v)&=&:\exp\left(
\sum_{m\neq 0}\frac{1}{m}(-B_m^j+B_m^{j+1})
(q^{N-j}z)^{-m}
\right):
e^{-\hat{\alpha}_j}z^{-h_j}(q^{-j+N}z)^{\frac{P_{\alpha_j}-1}{r}
+\frac{h_j}{r}},\label{free:F}
\end{eqnarray}
together with $H_j^\pm(v), K_j(v)$ given in
(\ref{def:total3c})-(\ref{def:total4c}) satisfy
the commutation relations in Definition
\ref{uqpslnh} for level $c=1$. 
\end{prop}

Now using 
this free field realization in \eqref{kjp}-
\eqref{elj}, we
obtain a realization of 
the $L$-operator $\hL^+(v)$ for $c=1$.
Using this in the "intertwining relation"
\eqref{Inter:Uqp2}, 
 we can solve it for the
II vertex operator.
\begin{thm}\lb{vertexop}
The highest components of the 
type II vertex operator 
$\Psi_N^*(z)$ is realized in terms of a free field by
\begin{eqnarray}
\Psi_N^*(z)&=&:\exp\left(
\sum_{m \neq 0}\frac{[r m]}{m[r^* m]}B_m^N z^{-m}
\right):
e^{-\bar{\Lambda}_{N-1}}z^{-h_{\bar{\epsilon}_N}}
e^{Q_{\bar{\epsilon}_N}}
z^{\frac{1}{r^*}P_{\bar{\epsilon}_N}}
z^{(1+\frac{1}{r^*})\frac{N-1}{2N}}
,\nonumber\\
\label{free:Psi}
\end{eqnarray}
where $\bar{\Lambda}_{N-1}=\frac{1}{N}
(\hat{\alpha}_1+2\hat{\alpha}_2+
\cdots +(N-1)\hat{\alpha}_{N-1}).$
The other components of the 
type II vertex $\Psi_j^*(z)\ (1\leq j \leq N)$ are given by
\begin{eqnarray}
\Psi_j^*(z)&=&
a_{j,N}^* \oint_{C^*} 
\prod_{m=j}^{N-1}
\frac{dw_m}{2\pi i w_m}
\Psi_N^*(v)
E_{N-1}(v_{N-1})\cdots
E_{j}(v_{j})
\nonumber\\
&&\times
\prod_{m=j}^{N-1}
\frac{[v_{m+1}-v_{m}-
P_{j,m+1}
+\frac{1}{2}]^*[1]^*}
{[v_{m+1}-v_{m}-\frac{1}{2}]^*
[P_{j,m+1}-1]^*}.\label{Type-II}
\end{eqnarray}
Here $v_N=v$. 
The integration contour $C^*$ is specified as follows.
For the integration contour for $w_{m}\ (j\leq m\leq N-1)$, the poles 
at $w_m=p^{*n}q^{-1}w_{m+1}\ (n=0,1,..)$ are inside, 
whereas  the poles at $w_m=p^{*-n} q w_{m+1}\ (n=0,1,...)$ are outside.
\end{thm}
\begin{thm}\label{thm:Vertexcom}~~~~
The free field realizations of 
the type-II vertex operator $\Psi_{\mu}^*(z)$
satisfies the following commutation relation. 
\begin{eqnarray}
&&{\Psi}^{*}_{j_1}(z_1)
{\Psi}^{*}_{j_2}(z_2)
=\sum_{j_1',j_2'=1}^N
{\Psi}^{*}_{j_2'}(z_2)
{\Psi}^{*}_{j_1'}(z_1)\ 
{R}^{*j_1j_2}_{j_1'j_2'}(u_1-u_2,P)
\label{Com:Type-II}
\end{eqnarray}
Here we set 
${R}^*(v,P)={\mu^*(v)}\bar{R}^{*}(v,P)$ 
with
\begin{eqnarray}
&&\mu^*(v)=z^{({\frac{1}{r^*}-1)}\frac{N-1}{N}}
\frac{\{pq^{2N-2}z\}^*
\{q^2z\}^*
\{p/z\}^*
\{q^{2N}/z\}^*
}{
\{pz\}^*
\{q^{2N}z\}^*
\{pq^{2N-2}/z\}^*
\{q^2/z\}^*
}.
\end{eqnarray}
\end{thm}

\subsection{The dual vertex operator}
The dual of the type II vertex operator of $U_{q,p}(\slnh)$ 
is an operator satisfying 
\bea
&&\Psi(z)\ :\widehat{\F} \to V_z\otimes \widehat{\F}'.
\ena
We define its components in the vector representation by 
\bea
&&\Psi(z)=\sum_{j=1}^N \bfv_j\otimes \Psi_j(z).
\ena
The following inversion relations hold.
\bea
&&\Psi_j(z)\Psi^*_k(z')=\delta_{j,k}\frac{g_N\ z^{1-N} }{1-q^{-N}\frac{z'}{z}}+\cdots,\\
&&g_N=\sqrt{-1}^{N}q^{\frac{N+1}{2r^*}+\frac{N^2-1}{2}}\left(\frac{(p^*q^2;p^*)_\infty}{(p^*;p^*)_\infty}\right)^{N}\frac{(pq^{2N};q^{2N},p^*)_\infty(q^{2N}q^{-2};q^{2N},p^*)_\infty}{
(q^{2N}p^{*};q^{2N},p^*)_\infty(q^{2N};q^{2N},p^*)_\infty},\nn
\ena
as $z'\to zq^{{N}}$, as  well as 
\bea
&&\sum_{j=1}^N\Psi_j(z)\Psi_j^*(z')=\frac{g'_N\ z^{1-N}}{1-q^N\frac{z'}{z}}+\cdots,\qquad \sum_{j=1}^N\Psi^*_j(z)\Psi_j(z')=\frac{g'_N\ z^{1-N}}{1-q^N\frac{z'}{z}}+\cdots,
\ena
where
\be 
&&g'_N=\sqrt{-1}^{-N}\frac{q^{-\frac{N+1}{2r^*}-\frac{N^2-1}{2}}}{(p^*;p^*)_\infty^{2N-3}
(q^{-2};p^*)_\infty^N}\frac{(p;q^{2N},p^*)_\infty(q^{-2};q^{2N},p^*)_\infty}{
(p^*;q^{2N},p^*)_\infty(q^{2N};q^{2N},p^*)_\infty},
\en
as $z' \to zq^{-{N}}$.
The free field realizaton is given  
as follows.
\begin{eqnarray}
\Psi_j(z)&=& 
\oint_{C} \prod_{m=1}^{j-1}
\frac{dz_m}{2\pi i z_m}
\Psi_1(z)E_{1}(v_{1})
\cdots E_{j-1}(v_{j-1})\nn\\
&&\qquad\times
\prod_{m=1}^{j-1}
\frac{[v_{m-1}-v_{m}-
P_{m-1,j}
+\frac{1}{2}]^*[1]^*}{[v_{m-1}-v_{m}-\frac{1}{2}]^*
[P_{m-1,j}-1]^*},
\qquad (1\leq j \leq N)\label{Type-I}
\ena
where $v=v_0$ and 
\bea
\Psi_1(z)&=&:\exp\left(-\sum_{m\not=0}\frac{[rm]}{m[r^*m]}B_m^1(q^Nz)^{-m}
\right):\ e^{\bar{\Lambda}_1}z^{h_{\bep_1}}e^{-Q_{\bep_1}}(q^Nz)^{-\frac{1}{r^*}
P_{\bep_1}+\frac{N-1}{2Nr^*}}z^
{\frac{N-1}{2N}},\nn\\
\end{eqnarray}
with $\bar{\Lambda}_1=\frac{1}{N}((N-1)\hat{\alpha}_1+(N-2)\hat{\alpha}_2+\cdots+
\hat{\alpha}_{N-1})$.
The integration contour $C$ 
is specified by the condition : for the contour for $w_m\ (1\leq m\leq j-1)$, the poles at $w_m=q^{-1}w_{m-1}p^{*n}\ (n=0,1,2,...)$ are inside, whereas the poles at $w_m=qw_{m-1}p^{*-n}\ (n=0,1,2,...)$ are outside. 

\noindent
{\it Remark}~~~The free field realizations of the vertex operators in Theorem
\ref{vertexop} and of the dual vertex operators 
are essentially the same as those of the $\slnh$ RSOS model 
obtained in \cite{AJMP,FKQ}. 

\section{Fusion of the Vertex Operators}
We now consider the fusion of the type II vertex operator $\Psi^*_1(z_2)$ and 
its dual $\Psi_1(z_1)$. Namely, we consider a product $\Psi_1(z_1)\Psi_1^*(z_2)$ and investigate the limits to the fusion points $z_1= q^{-N}p^{*n}z_2\ 
(n=0,1,2,..,N)$, where the contour in \eqref{Type-II} for $w_1$ gets pinches.

For example, let us consider the case $n=1$. If we take residues for the poles 
$w_{N-1}=q^{-1}z_2,\ w_{j-1}=q^{-1}w_{j}\ (j=N-1,N-2,..,3)$, the limit $z_1\to 
q^{-N}p^*z_2$ causes pinches in the contour for $w_1$ at two points $w_1=q^{-(N-1)}z_2,\ q^{-(N-1)}p^*z_2$. Similarly, for $1\leq l\leq N-2$, 
if we take residues at the poles
$w_{N-1}=q^{-1}z_2,\ w_{j-1}=q^{-1}w_{j}\ (j=N-1,N-2,..,N-l+1),\ w_{N-l}=
q^{-1}p^*w_{N-l+1},\ w_{j-1}=q^{-1}w_{j}\ (j=N-l-1,N-l-2,..,3)$, 
the same limit $z_1\to q^{-N}p^*z_2$ causes a pinch in the contour for $w_1$ at a point $w_1=q^{-(N-1)}p^*z_2$. Hence in the limit $z_1\to q^{-N}p^*z_2$, 
we obtain totally $N$ terms of contributions from the residues at the $N$ 
pinching points. 
Similar consideration leads us to the following results. As  
${z_1\to q^{-N}p^{*n}z_2}$, 
\bea
&&\Psi_1(z_1)\Psi_1^*(z_2)\nn\\
&&=\frac{z_1^{1-N}}{1-
q^{-N}p^{*n}\frac{z_2}{z_1}}\left\{C_n\ \tilde{T}_n(q^{(n-1)r^*}z_2)\right.\nn\\
&&+\left.{\sum}'_{1\leq j_1\leq j_2\leq ..
\leq j_n\leq N}C_{j_1,j_2,..,j_n}:\Lambda_{j_1}(z_2q^{(2n-1)r^*})\Lambda_{j_1}(z_2q^{(2n-3)r^*})\cdots \Lambda_{j_n}(z_2q^{r^*}):\right\}+\cdots.\nn\\
&&\lb{fusion}
\ena
Here 
\bea
&&\tilde{T}_n(z)={\sum}_{1\leq j_1< j_2< ..
< j_n\leq N}:\Lambda_{j_1}(z_2q^{(n-1)r^*})\Lambda_{j_2}(z_2q^{(n-3)r^*})\cdots \Lambda_{j_n}(z_2q^{-(n-1)r^*}):,\\
&&\Lambda_j(z)=:\exp\left(\sum_{m\not=0}\frac{q^{rm}-q^{-rm}}{m}B^j_mz^{-m}
\right): q^{-2P_{\bep_j}}p^{*h_{\bep_j}}q^{\frac{2(1-N)}{N}}p^{*-\frac{1}{N}-j},\\
&&C_n=\sqrt{-1}^Nq^{\frac{N+1}{2r^*}+\frac{N^2-1}{2}}\left(
\frac{(p^*q^2;p^*)_\infty}{(p^*;p^*)_\infty}\right)^N\left(
\frac{1-pq^{-N}}{1-q^{-N}}\right)^n\nn\\
&&\qquad\qquad \times\frac{(pq^{2N}p^{*-n};q^{2N},p^*)_\infty(q^{2N-2}p^{*-n};q^{2N},p^*)_\infty}
{(q^{2N}p^{*-n};q^{2N},p^*)_\infty(q^{2N}p^{*1-n};q^{2N},p^*)_\infty}.\nn\\
\ena   
In \eqref{fusion}, ${\sum}'$ denotes the sum over the complementary set to $1\leq j_1< j_2< ..
< j_n\leq N$. $C_{j_1,j_2,..,j_n}$ are constants not important here. 

The basic operators $\Lambda_j(z)\ (1\leq j \leq N-1)$ coincides with those in
 the  deformed $W_N$ algebra\cite{FF,AKOS}. 
The expressions for $\tilde{T}_n\ (1\leq n\leq N)$ are almost 
 same as 
those of the generating ``currents" of the  deformed $W_N$ algebra, but 
 the unit of the $q$-shift in the arguments in $\Lambda_j(z)$ is different. 
 In an 
identification of the parameters $p_W=q^{-2},\ q_W=p=q^{2r}$, where $p_W$ 
and $q_W$ are $p$ and $q$ in \cite{FF,AKOS}, respectively, the unit of the $q$-shift in \cite{FF,AKOS}
is given by $p_W$, whereas it is $p^*=q^{2(r-1)}$ in our $\tilde{T}_n(z)$. 
As a consequence, we have 
\bea
&&\tilde{T}_N(z)=:\Lambda_{1}(z_2q^{(N-1)r^*})\Lambda_{2}(z_2q^{(N-3)r^*})\cdots \Lambda_{N}(z_2q^{-(N-1)r^*}):\not=1.
 \ena  
Therefore, our deformed $W$ algebra generated by $\tilde{T}_n\ (1\leq n\leq N)$ is  
$\goth{gl}_N$ type instead of $\sln$ type. 

On the other hand, since the type II vertex operator 
$\Psi^*(z)$ and its dual $\Psi(z)$ are the creation operators of the physical excited particle and anti-particle, it is natural to identify the operators  
$\tilde{T}_n(z)\ (1\leq n \leq N)$ with the creation operator of 
their bound states. The $S$-matrix of the bound state particles are calculated as follows. 
\bea
&&\tilde{T}_n(z)\tilde{T}_m(w)=S_{n,m}(w/z)\ \tilde{T}_m(w)\tilde{T}_n(z),\\
&&S_{n,m}(z)=\prod_{k=1}^n\prod_{l=1}^m \varphi_N\left(
{z}q^{r^*(n-m+2(l-k))}\right),\\
&&\varphi_N(z)=\frac{\Theta_{q^{2N}}(q^{2}z)\Theta_{q^{2N}}(p^{*}z)\Theta_{q^{2N}}(p^{*-1}q^{-2}z)}{\Theta_{q^{2N}}(q^{-2}z)\Theta_{q^{2N}}(p^{*-1}z)\Theta_{q^{2N}}(p^{*}q^{2}z)}.
\ena
Again, this $S$-matrix is different from the one obtained by Feigin and Frenkel (sec7.2 in \cite{FF} ) only by the choice of the unit of the $q$-shift.

The scaling limit of the $\slnh$ RSOS model is expected to be the RSOS 
restriction of the affine Toda field theory with imaginary coupling constant.
It is interesting to compare the scaling 
limit of our $S$-matrices, $R^*(v,P)$ for the excited particle and 
$S_{n,m}(z)$ for 
the bound states, with the bootstrap results 
\cite{Johnson,Ganden}.



~\\
{\bf Acknowledgements}~~
We would like to thank Patrick Dorey, Tetsuji Miwa and Robert Weston 
for discussion. We also thank the organizers of RAQIS'03, 
Daniel Arnaudon, Jean Avan, Luc Frappat, Eric Ragoucy and Paul Sorba, for 
their kind invitation to the conference and their hospitality.  
H.K. is also  grateful to DAMTP, 
University of Cambridge, for hospitality. 
This work is  supported by the JSPS/Royal Society fellowship
and the Grant-in-Aid for Young Scientist ({\bf B}) (14740107)
 from the Ministry of Education, Japan .


\end{document}